# *Master Stability Function for networks of coupled non-smooth oscillators*

Volodymyr Denysenko[1,a)], Marek Balcerzak[1], Artur Dabrowski[1]

## *Affiliations*

[1] Division of Dynamics, Lodz University of Technology, 90-924 Lodz, Poland

[a] Author to whom correspondence should be addressed: volodymyr.denysenko@dokt.p.lodz.pl

## Abstract

This paper describes the extension of the Master Stability Function(MSF) for arrays of non-smooth oscillators. This extension is based on the previously introduced Jacobi matrix estimation method. The extended method can be applied to networks of non-smooth coupled oscillators. The extension and its limitations are described. Then, a new algorithm is applied to calculate the MSF for an impact oscillator. The results are presented and validated using an alternative MSF estimation method.

Keywords: Master Stability Function, non-smooth, complex networks, nonlinear dynamics, stability, complete synchronization.

**Complete synchronization (or the lack thereof) in arrays of oscillators plays a significant role in various fields of science and engineering [1] [2] [3] [4] [5]. To assess the stability of complete synchronization in oscillator arrays while simplifying the problem to a low-dimensional variational equation analysis, a method called the Master Stability Function (MSF) was introduced. This method has been successfully applied to networks of smooth oscillators. However, to ensure accuracy of a model, many processes and phenomena across various scientific fields must be represented using non-smooth differential equations [6] [7] [8] [9]. Yet, applying the MSF to non-smooth oscillator arrays is more challenging and requires an alternative estimation technique. This article presents a novel approach to MSF calculation that can be applied to the analysis of non-smooth oscillator arrays.**

## Introduction

The investigation of collective behavior in arrays of oscillators has gathered significant attention due to its relevance in various fields of science and engineering [1] [2] [3] [4] [5]. Oscillatory systems, when interconnected, often exhibit complex dynamic behaviors that emerge from their interactions. One of the phenomena that was observed in such networks is synchronization [10], where all oscillators in the array oscillate with the same frequency. Aforementioned coupled system behavior is not only a fundamental phenomenon in nonlinear dynamics theory, but has also been observed in practice [11] [12] [13] [14] [15]. The scope of this article is limited to complete synchronization – identical motion of oscillators [16] [17].

The Master Stability Function (MSF) [18] is one of the most commonly applied methods for investigating the stability of complete synchronization in networks of identical oscillators. This approach relies on analyzing a low-dimensional variational equation around the synchronization manifold. The MSF has two significant advantages. First, results of computations for one particular choice of oscillator and selected coupling configuration can be reused for various possible connection topologies. Second, its similarity to another significant and widely applied trajectory stability criterion, the Largest Lyapunov Exponent (LLE) [19], allows MSF to be estimated using well-known methods for calculating the LLE [20] [21] [22] [23] [24]. These advantages make the MSF a powerful tool for investigating networks of continuous oscillators. However, the estimation of the MSF for a network of non-smooth oscillators is a much more complex task because the Jacobian matrix cannot be calculated at non-differentiable points, making the variational equation undefined. Nevertheless, many different systems and processes in various fields of science and engineering, including mechanical oscillators with dry friction [6], electro-mechanical systems [7], brain functioning [8] [9], can be accurately described using non-smooth models. This implies the need to extend the MSF in a way that enables the investigation of such systems.

To overcome difficulties related to MSF calculation for arrays of non-smooth oscillators, several approaches have been proposed. Some of these methods, such as the two-oscillators [25] and the three-oscillators probe [26], allow the determination of the sign of MSF, which is sufficient for assessing the stability of complete synchronization. These methods can be applied to both continuous and non-smooth oscillators, although they may require longer computation times. Additionally, all possible interactions near non-differentiable points between oscillators must be considered, which can complicate the calculation process. A different approach, which involves direct calculation of the MSF, was introduced in [27], where it was suggested to calculate and use saltation matrix for mapping perturbed trajectories through a switching manifold. Another method, based on calculating the monodromy matrix for MSF estimation was proposed in [28]; however, this approach is suitable only for Filippov systems and requires existence of limit cycle for single oscillator.

This paper introduces a novel technique for calculating the MSF for non-smooth oscillators, inspired by a method used for Lyapunov Exponent (LE) estimation in non-smooth systems [29]. The core concept of the aforementioned algorithm is the calculation of the Jacobi matrix using initial perturbations in the form of orthogonal vectors. The modification introduced in this paper enables the estimation of the Jacobi matrix around non-differentiable points for piecewise-smooth systems, thereby allowing the calculation of the MSF for the network of oscillators under investigation.

## *Method*

To introduce the MSF approach, let us consider a network system of N identical oscillators, which can be described by Eq. (1), where $\boldsymbol{x} = (x_1, x_2, \dots, x_N)$ is a vector representing the state of each oscillator in the array, $\boldsymbol{F}(\boldsymbol{x}) = [F(x_1), F(x_2), \dots, F(x_n)]$ is a vector of vector fields that defines the dynamics of each individual oscillator in the phase space, $\boldsymbol{G}$ is the matrix of coupling coefficients, $\boldsymbol{H}(\boldsymbol{x}) = [H(x_1), H(x_2), \dots, H(x_n)]$ is a vector of functions of each node's variable used in the coupling and $\otimes$ represents the Kronecker product.

$$\dot{\boldsymbol{x}} = \boldsymbol{F}(\boldsymbol{x}) + \sigma \boldsymbol{G} \otimes \boldsymbol{H}(\boldsymbol{x}) \qquad (1)$$

Then, by obtaining the variational equation of Eq. (1) and diagonalizing the matrix **G**, a block-diagonalized variational equation is derived, with each block having the following form:

$$\dot{\boldsymbol{\xi}}_k = [D\boldsymbol{F} + \sigma \gamma_k D\boldsymbol{H}] \boldsymbol{\xi}_k \qquad (2)$$

where $\gamma_k$ is an eigenvalue of **G**, and $\boldsymbol{\xi}_k$ is the eigenvector corresponding to the $k$-th eigenvalue, $k \in \{0, 1, \dots, N-1\}$. It is worth mentioning that $\gamma_0 = 0$, therefore for this eigenvalue Eq. (2) reduces to the variational equation of a single, separated node. The other eigenvalues correspond to transverse eigenvectors. Substituting $\sigma\gamma = \alpha + i\beta$, where $\gamma$ represents any arbitrary eigenvalue of $G$, the variational equation is obtained.

$$\dot{\xi} = [DF + (\alpha + i\beta)DH]\xi \qquad (3)$$

In order to evaluate the stability of the complete synchronization in the investigated network, the largest Lyapunov exponent (LLE) for Eq. (3) must be calculated as a function of α and β. In this case, LLE indicates the divergence of nearby trajectories in directions transversal to the synchronization manifold; therefore, it is also commonly referred to as Transversal Lyapunov Exponent (TLE) [30]. Then, for any given σ, a corresponding value of TLE can be assigned for each $\gamma_k$, and stability of each eigenmode can be determined. For a stable synchronization state of the oscillators network, each eigenmode must correspond to a negative Lyapunov exponent.

As demonstrated above, the estimation of the MSF is based on a variational equation similar to that one on which the LLE calculation is based. Therefore, in the case of a network of smooth oscillators, the MSF can be estimated using classical methods for the LLE estimation. However, in the case of a non-smooth oscillator, the expression inside the parentheses of Eq. (3) cannot be calculated at points of non-differentiability, making it impossible to directly predict the evolution of a perturbation vector. To overcome this problem in case of LLE calculations, several methods have been proposed [31] [32] [33] [34] [35] [36] [37] [38] [39]. However, we will focus on the method that allows us to calculate the Jacobi matrix of the trajectory [29].

To briefly introduce this method, let's consider an autonomous n-dimensional dynamical system, described by equation:

$$\dot{x}(t) = f(x(t)) \qquad (4)$$

where $x \in \mathbb{R}^n$ is the system's state and $f \in \mathbb{R}^n$ is the vector field describing the dynamics of the system. Let $\phi_t(x_0)$ be the solution of the Eq. (4) starting from the initial conditions $x_0$. Then, the Jacobi matrix is as follows.

$$\Phi_t(x_0) = \frac{d\phi_t(x_0)}{dx_0} \qquad (5)$$

The most important property of the Jacobi matrix $\Phi_t(x_0)$ in the context of calculating the Lyapunov exponent is that it describes the evolution of an infinitesimal perturbation of the initial conditions $\delta x_0$ over a time $t$ [22].

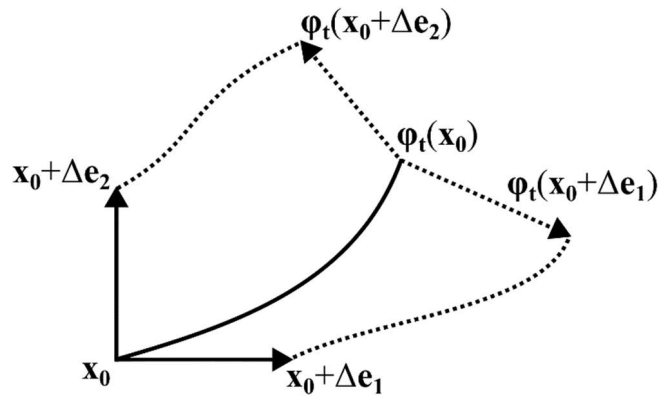

*Figure 1 Illustration of the algorithm of the $\Phi(x)$ matrix approximation for $f: \mathbb{R}^2 \to \mathbb{R}^2$*

It is proposed to interpret $\boldsymbol{\Phi}_t$ as the Jacobi matrix of the trajectory $\boldsymbol{\phi}_t$ evaluated at the point $\boldsymbol{x_0}$, which allows it to be represented as:

$$\boldsymbol{\Phi}_t(\boldsymbol{x_0}) = \frac{d\boldsymbol{\phi}_t(\boldsymbol{x_0})}{d\boldsymbol{x_0}} = \lim_{\Delta \to 0} [\frac{\boldsymbol{\phi}_t(\boldsymbol{x_0} + \Delta \boldsymbol{e_1}) - \boldsymbol{\phi}_t(\boldsymbol{x_0})}{\Delta}, \ldots, \frac{\boldsymbol{\phi}_t(\boldsymbol{x_0} + \Delta \boldsymbol{e_n}) - \boldsymbol{\phi}_t(\boldsymbol{x_0})}{\Delta}] \qquad (6)$$

where $\boldsymbol{e_1}, \ldots, \boldsymbol{e_n}$ are vectors that form an orthonormal basis in $\mathbb{R}^n$. To estimate the Jacobi matrix numerically for some initial conditions $\boldsymbol{x_0} \in \mathbb{R}^n$, it is proposed to perturb the initial conditions by a sufficiently small value $\Delta$ in the direction of each basis vector of $\mathbb{R}^n$, which results in a set of initial conditions: $\boldsymbol{x_0}, \boldsymbol{x_0} + \Delta \boldsymbol{e_1}, \ldots, \boldsymbol{x_0} + \Delta \boldsymbol{e_n}$. Next, the solution for each initial condition at the same moment of time $t$ must be evaluated. Finally, using the obtained trajectories, the Jacobi matrix can be approximated as follows.

$$\boldsymbol{\Phi}_t(\boldsymbol{x_0}) \approx \frac{1}{\Delta}[\boldsymbol{\phi}_t(\boldsymbol{x_0} + \Delta\boldsymbol{e_1}) - \boldsymbol{\phi}_t(\boldsymbol{x_0}), \ldots, \boldsymbol{\phi}_t(\boldsymbol{x_0} + \Delta\boldsymbol{e_n}) - \boldsymbol{\phi}_t(\boldsymbol{x_0})] \qquad (7)$$

A graphical representation of the procedure is shown in Figure 1. This approach allows for the calculation of the Jacobi matrix for any period of time for both smooth and non-smooth systems. However, in the case of non-smooth systems, it is necessary that there is no discontinuity between the $\boldsymbol{\phi}_t(\boldsymbol{x_0})$ and $\boldsymbol{\phi}_t(\boldsymbol{x_0} + \Delta \boldsymbol{e_i})$ for $i = 1, \ldots, n$. The calculated Jacobi matrix can be applied to estimate the perturbation behavior over the chosen time period. Unfortunately, this method cannot be used directly to calculate the MSF; therefore, we propose modifications to adapt this algorithm for this purpose.

In order to develop an approach that allows adjusting the Jacobi matrix of the trajectory of a single oscillator to calculate the evolution of perturbation, described by Eq. (3), we need to assume that the evolution of the perturbation can be linearized around the non-differentiable point in phase space. Moreover, it is required that the Jacobi matrix of the trajectory $\boldsymbol{\Phi}_t(\boldsymbol{x_0})$ can be expressed in the form $\exp(At)$, i.e., the logarithm $A$ of the matrix $\boldsymbol{\Phi}_t(\boldsymbol{x_0})$ is well-defined.

$$\delta x(t) = \boldsymbol{\Phi}_t(\boldsymbol{x_0})\delta \boldsymbol{x_0} = e^{At}\delta \boldsymbol{x_0} \qquad (8)$$

However, the Jacobi matrix, corresponding to the solution of Eq. (3) is given by:

$$\boldsymbol{\Phi}_t(\boldsymbol{x_0}) = \exp([A + (\alpha + i\beta)H]t) \qquad (9)$$

The problem to be solved assumes that $\boldsymbol{\Phi}_{At}(\boldsymbol{x_0})$ for an isolated oscillator and $(\alpha + i\beta)H$ are known, and $\boldsymbol{\Phi}_t(\boldsymbol{x_0})$ for a coupled oscillator is to be determined. The algorithm for solving this problem is as follows: the matrix $At$ is calculated as the logarithm of $\boldsymbol{\Phi}_{At}(\boldsymbol{x_0})$, and then $[A + (\alpha + i\beta)H]t$ must be estimated, which is sufficient to determine $\boldsymbol{\Phi}_t(\boldsymbol{x_0})$.

However, it is worth to mention a few important aspects of the described method that should be considered. The first concerns the logarithm of matrix – it is important to note that only invertible matrices can have a logarithm [40]. In the context of the described calculations, this implies that the reduction in the order of the dynamical system does not occur. Therefore, the proposed method cannot be applied to systems analogous to the stick-slip oscillator [41], where the system's motion is governed by two different vector fields of different orders on either side of the discontinuous manifold [42].

The second and final aspect to consider is ensuring that the evolution of perturbations can be linearized around non-differentiable points. To achieve this, the calculation time for the Jacobi matrix near a discontinuity must be minimized—ideally, it should be on the order of the integration step size. When this condition is met, the technique described above can be applied to estimate the Jacobi matrix that characterizes the evolution of perturbations, as described by Eq. (3), across the non-differentiable points. This enables the estimation of the TLE, and consequently, the MSF for non-smooth oscillators.

## Numerical example

The proposed procedure for calculating the MSF for non-smooth oscillators was applied to a simple mechanical oscillator with impacts. The system is a typical linear mass-spring-damper system, which can collide with a wall placed at the position $X_w$. The equation of motion for the system is as follows:

$$\dot{X}_1(t) = X_2(t)$$
$$m\dot{X}_2(t) = -cX_2(t) - kX_1(t) + F * \cos(\Omega t)$$
$$\dot{X}(t_c^+) = -R\dot{X}(t_c^-), X(t_c) = X_w \quad (9)$$

where $c$ is the damping coefficient, $k$ is the spring stiffness, $F$ and $\Omega$ are the amplitude and frequency of the external forcing, respectively. $R$ is the coefficient of restitution, and $t_c$ is the moment in time when the collision with the wall occurs. Eq. (9) can be reduced to a dimensionless form:

$$\dot{x}_2(\tau) = x_2(\tau)$$
$$\dot{x}_2(\tau) = -2\zeta x_2(\tau) - x_1(\tau) + f\cos(\eta \tau)$$
$$\dot{x}(\tau_c^+) = -R\dot{x}(\tau_c^-), x(\tau_c) = x_w \quad (10)$$

where $\zeta = \frac{c}{2\sqrt{mk}}$ is the damping ratio, $\omega = \sqrt{\frac{k}{m}}$ is the natural frequency of oscillator, $f = \frac{F}{m\omega^2}$, $\tau = \omega t$ is the dimensionless time, $\eta = \frac{\Omega}{\omega}$ is the dimensionless frequency of external excitation, $x_w = \frac{kX_w}{F}$, and $\tau_c$ is the dimensionless time of collision. Two cases where considered: elastic and non-elastic collision. In both cases, the assumed coupling is described by the matrix:

$$\mathbf{H} = \begin{bmatrix} 0 & 0 \\ 1 & 0 \end{bmatrix} \quad (11)$$

which can be considered as a linear spring coupling, one of the most commonly encountered in mechanical systems. The coupling strength is considered as a variable parameter. Another assumption regarding the coupling is that the interaction between each pair of nodes is symmetrical; therefore, all eigenvalues of the coupling coefficient matrix are real. For the elastic collision, the system parameter values were chosen as follows: R = 1.0, $x_w$ = 2.0, $\zeta = 0.05, \eta = 0.712$. For an inelastic collision, the system parameter values were: R = 0.9, $x_w$ = 1.5, $\zeta = 0.05, \eta = 0.5975$. Since the system is piecewise linear, it was solved analytically. The occurrence of a collision was checked at regular intervals ($\Delta t = 10^{-3}$). To determine the collision time with proper accuracy, the bisection procedure was applied. To reduce transient effects, TLE calculations were started after a time interval equal to 500 periods of external excitation. Since the parameters of the single oscillator remained the same, this procedure was performed once, and then the oscillator state could be used for each parameter value. TLE calculations were carried out using the method described in the previous section.

It is worth mentioning that a real matrix can have a complex logarithm [40]. In such cases, the complex Jacobi matrix for MSF will be obtained, which will result in real perturbation vector becoming complex. However, since perturbation is expected to be real in a real system, the complex part of the calculated Jacobi matrix for TLE estimation was treated as a numerical error arising from linearization of the Jacobi matrix of the single system and discarded. TLE was calculated until one of the following conditions was satisfied: either the standard deviation of the last 100 calculated values of TLE was below a predefined threshold of $10^{-5}$, or the simulation time exceeded the limit of 2000 periods of external excitation.

In order to verify the results obtained with the proposed procedure, the two-oscillators probe was applied. The system of two oscillators was also solved analytically due to the linearity of the system. The procedure for collision detection and collision time determination was the same as in the case of TLE calculations. The initial length of perturbation vector for two-oscillator probe was set to $10^{-3}$; the direction of initial perturbation vector was chosen randomly. The termination conditions for the calculations were as follows: the occurrence of complete synchronization between the oscillators or the exceeding of the simulation time limit, which was also set to 2000 periods of external excitation. The local maxima of the difference in oscillator positions over time, taken from the last 100 periods of the simulation, were used to construct the bifurcation diagram. Figure 2 presents the bifurcation diagram and TLE graph for with an elastic collision, while Figure 3 illustrates the bifurcation diagram for the two-oscillator probe and the TLE graph for the scenario in which a system with an inelastic collision was analyzed.

## *Results and summary*

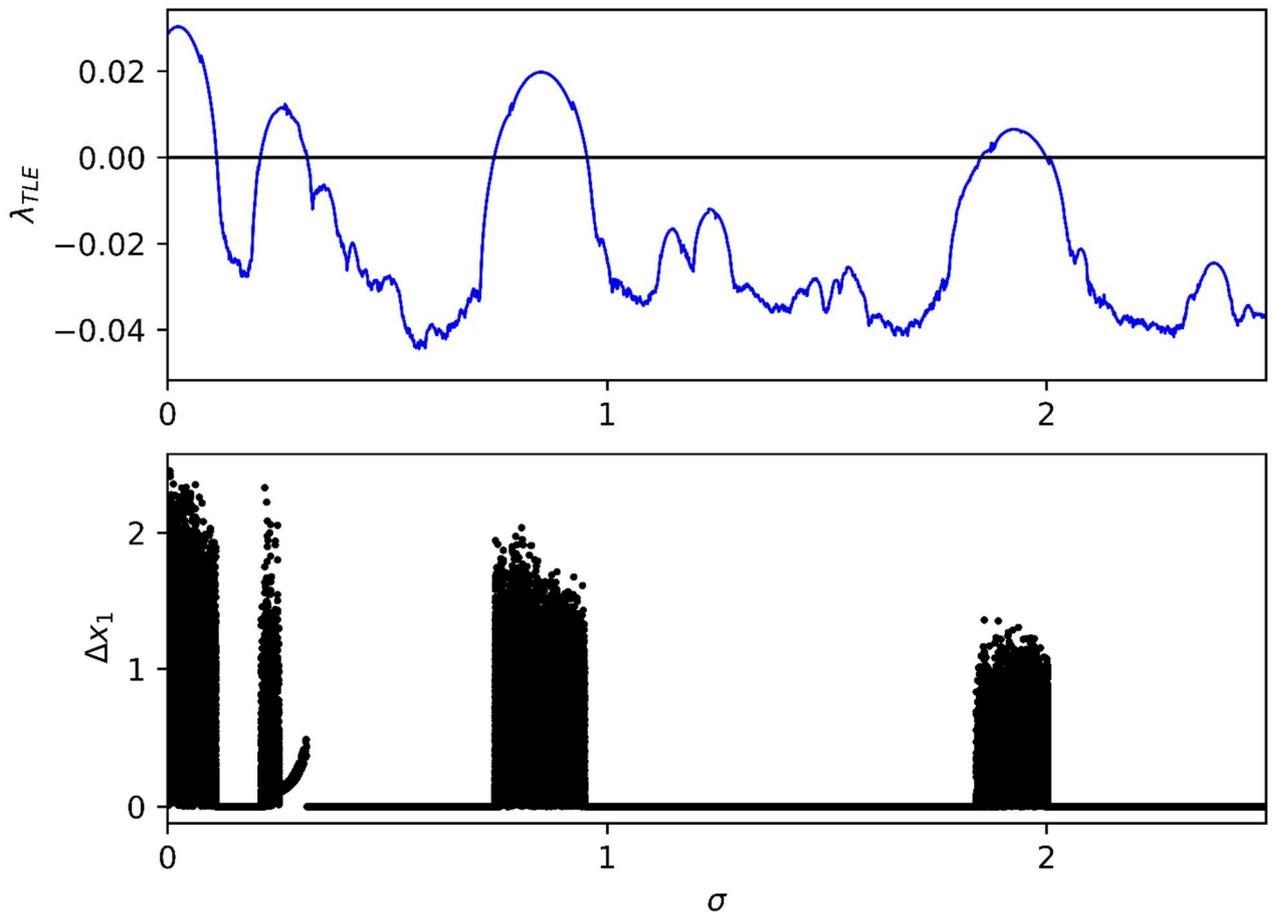

*Figure 2 TLE graph and bifurcation diagram for two-oscillators probe for oscillators with elastic collision*

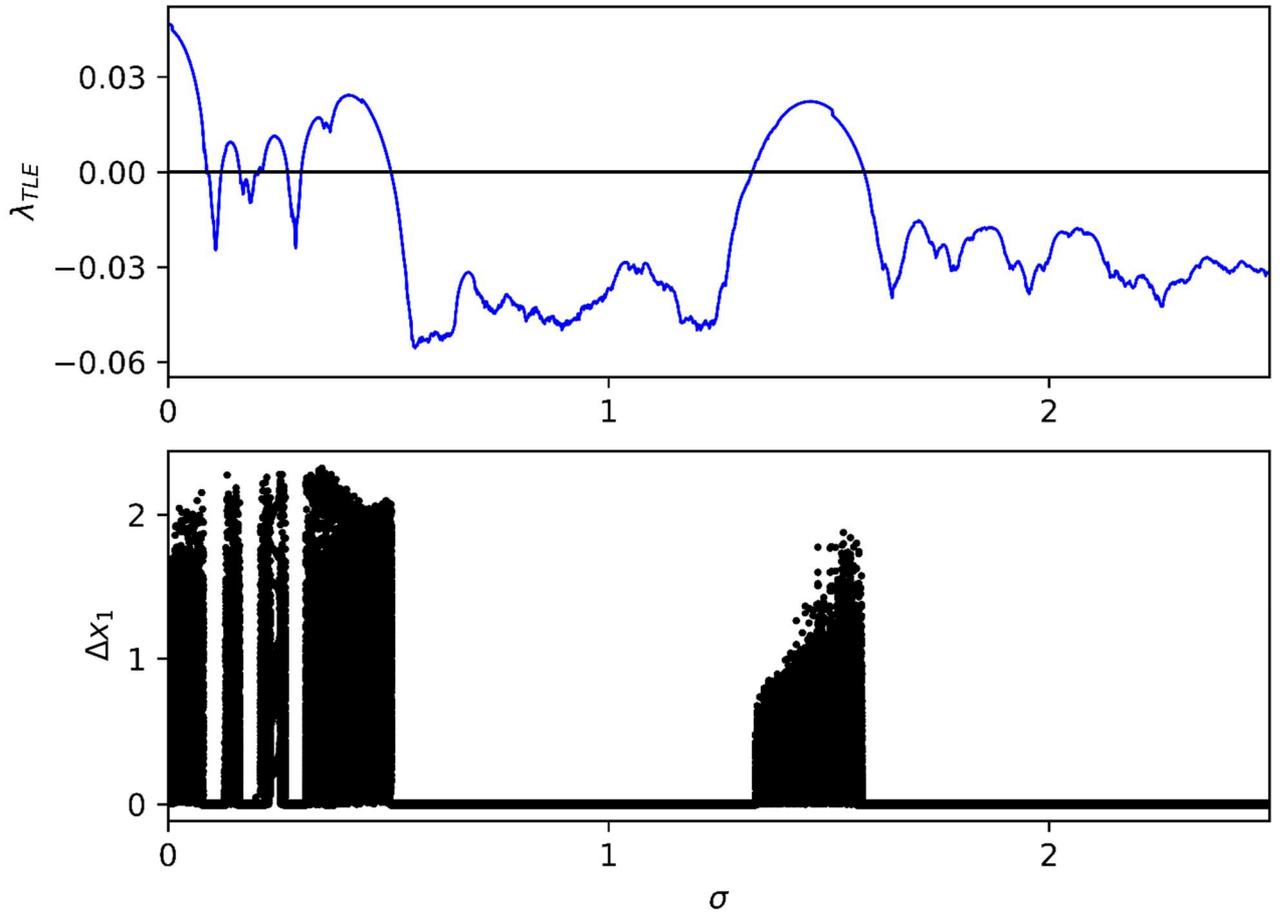

*Figure 3 TLE graph and bifurcation diagram for two-oscillators probe for oscillators with non-elastic collision*

Comparing the bifurcation diagrams obtained from the two-oscillator probe with the corresponding TLE graphs, we can conclude that the proposed method is sufficiently accurate, at least qualitatively. However, it is worth noting that the results obtained in regions where $\lambda_{TLE} \approx 0$ differ — complete synchronization occurs for positive values of $\lambda_{TLE}$ and vice versa. Nevertheless, such system behavior is expected.

This article introduces a new method for calculating the MSF in networks of non-smooth oscillators. The foundation of the method is presented, and the proposed calculation algorithm is outlined. Numerical examples are provided and compared with a well-known MSF calculation approach. Based on the presented examples, it can be concluded that the proposed method is applicable for MSF estimation in networks of non-smooth oscillators. The advantages of the proposed method are as follows: on one hand, it eliminates the need to account for the interactions of two or three oscillators near non-differentiable points in the calculations, which can be burdensome in non-smooth systems; on the other hand, it can be applied to almost any oscillator without requiring additional initial calculations.

## *Acknowledgements*

V.D has been supported by the National Science Centre, Poland, PRELUDIUM Program (Project No. 2023/49/N/ST8/02436) This paper has been completed while the first author was the Doctoral Candidate in the Interdisciplinary Doctoral School at the Lodz University of Technology, Poland.

## Author Declarations

## Conflict of interest

The authors have no conflicts to disclose

## Author contributions

**Volodymyr Denysenko:** Conceptualization; Investigation; Methodology(equal) Software; Visualization; Writing/Original Draft Preparation; Writing/Review & Editing (equal). **Marek Balcerzak:** Methodology(equal); Supervision(equal); Writing/Review & Editing (equal). **Artur Dabrowski:** Funding Acquisition; Supervision(equal)

## Data availability

The data that support the findings of this study are available from the corresponding author upon request.